\pdfoutput=1
\documentclass[reqno,12pt]{article}
\usepackage{amsmath}
\usepackage{amsthm}
\usepackage{amssymb}
\usepackage{authblk}
\usepackage{enumerate}
\usepackage{graphicx,tikz}
\usetikzlibrary{calc,math}
\newcommand{\subs}{\subseteq}
\newcommand{\setm}{\setminus}
\newcommand{\lf}{\lfloor}
\newcommand{\rf}{\rfloor}
\newcommand{\lc}{\lceil}
\newcommand{\rc}{\rceil}
\newtheorem{theorem}{Theorem} [section]
\newtheorem{lemma}[theorem]{Lemma}
\newtheorem{proposition}[theorem]{Proposition}

\newtheorem{cor}[theorem]{Corollary}
\newtheorem{conj}[theorem]{Conjecture}
\def\N{\mathbb N}
\def\R{\tilde R}
\def\RR{\text{RR}}
\def\rr{\tilde{r}}
\def\sr{\hat{r}}
\def\rn{\text{Red}}
\def\bn{\text{Blue}}
\def\ess{\text{ess}}
\newcommand{\ver}[2]{ 
\node [circle,inner sep=1pt,fill=black,draw,very thick] (#1) at (#2) {};
}

\newcommand{\edr}[2]{ 
\draw [red,very thick] (#1) to (#2);
}

\newcommand{\edb}[2]{ 
\draw [blue,very thick] (#1) to (#2);
}

\newcommand{\edar}[2]{ 
\draw [very thick,red] (#1) to [bend right=25] (#2);
}

\newcommand{\edab}[2]{ 
\draw [very thick,blue] (#1) to [bend right=27] (#2);
}

\newcommand{\bbbb}{ 
\foreach \X in {1,...,5}{\ver{a\X}{$(\X-1,0)$};}
\foreach \X in {2,...,5}{
\tikzmath{\Y=int(\X-1);}
\draw [blue,very thick] (a\X)--(a\Y);
}
}

\newcommand{\brrb}{ 
\foreach \X in {1,...,5}{\ver{a\X}{$(\X-1,0)$};}
\foreach \X in {2,...,5}{
\draw [blue,very thick] (a1)--(a2);
\draw [red,very thick] (a2)--(a3);
\draw [red,very thick] (a3)--(a4);
\draw [blue,very thick] (a4)--(a5);
}
}

\newcommand{\bbrb}{ 
\foreach \X in {1,...,5}{\ver{a\X}{$(\X-1,0)$};}
\foreach \X in {2,...,5}{
\draw [blue,very thick] (a1)--(a2);
\draw [blue,very thick] (a2)--(a3);
\draw [red,very thick] (a3)--(a4);
\draw [blue,very thick] (a4)--(a5);
}
}

\newcommand{\bbrr}{ 
\foreach \X in {1,...,5}{\ver{a\X}{$(\X-1,0)$};}
\foreach \X in {2,...,5}{
\draw [blue,very thick] (a1)--(a2);
\draw [blue,very thick] (a2)--(a3);
\draw [red,very thick] (a3)--(a4);
\draw [red,very thick] (a4)--(a5);
}
}

\newcommand{\rrbr}{ 
\foreach \X in {1,...,5}{\ver{a\X}{$(\X-1,0)$};}
\foreach \X in {2,...,5}{
\draw [red,very thick] (a1)--(a2);
\draw [red,very thick] (a2)--(a3);
\draw [blue,very thick] (a3)--(a4);
\draw [red,very thick] (a4)--(a5);
}
}

\newcommand{\limb}{ 
\ver{l1}{$(0,-1.5)$};
\ver{l2}{$(1,-1.5)$};
\ver{l3}{$(2,-1.5)$};
\ver{l4}{$(3,-1.5)$};
\ver{l5}{$(1,-2.5)$};
\edb{l1}{l2};
\edr{l2}{l3};
\edb{l3}{l4};
\edb{l2}{l5};
}
\begin{document}

\title{Off-diagonal online size Ramsey numbers for paths}
\author{Ma\l gorzata Bednarska-Bzd\c ega}
\affil{Faculty of Mathematics and CS, Adam Mickiewicz University, Pozna\'n, Poland}

\maketitle

\begin{abstract}
Consider the following Ramsey game 
played on the edge set of $K_\N$. In every round, Builder selects an edge and Painter colours 
it red or blue.  Builder's goal is to force Painter to create a red copy of a path $P_k$ on $k$ vertices 
or a blue copy of $P_n$ as soon as possible. The online (size) Ramsey number $\rr(P_k,P_n)$ 
is the number of rounds in the game provided  Builder and Painter play optimally. 
We prove that $\rr(P_k,P_n)\le (5/3+o(1))n$ provided $k=o(n)$ and $n\to \infty$.
We also show that $\rr(P_4,P_n)\le\lc 7n/5\rc -1$ for $n\ge 10$, which improves the upper bound 
obtained by J.~Cyman, T.~Dzido, J.~Lapinskas, and A.~Lo
and implies their conjecture that $\rr(P_4,P_n)=\lc 7n/5\rc -1$. 
\end{abstract}

\section{Introduction}

Let $G$ and $H$ be finite graphs. Consider the following game  $\R(G,H)$ played on the infinite board $K_\N$ (i.e. the board is a complete graph with the vertex set $\N$). In every round, Builder chooses a previously unselected edge of $K_\N$ and Painter colours it red or blue. The game ends when there is a red copy of $G$ or a blue copy of a $H$ on the board.  Builder aims to finish the game as soon as possible, while Painter tries to avoid a red $G$ and a blue $H$ as long as possible. 
By $\rr(G,H)$ we denote the number of rounds in the game $\R(G, H)$, provided both players play optimally and we call it the
online size Ramsey number for $G$ and $H$. In the literature online size Ramsey numbers are called also online Ramsey numbers. The online size Ramsey numbers  $\rr(G, H)$ are game counterparts of well known size Ramsey numbers; the  size Ramsey number $\sr(G,H)$ is the minimum number of edges in a graph with the property that every two-colouring of its edges results in a red copy of $G$ or a blue copy of $H$. Clearly $\rr(G,H)\le \sr(G,H)$. 

In this paper we study online size Ramsey numbers for $G=P_k$ and $H=P_n$, where $P_t$ denotes a path on $t$ vertices.
Games in which Builder tries to force Painter to create a monochromatic path were considered also in other variants: the induced version of the online Ramsey number was studied in \cite{ind}, ordered path games on the infinite (ordered) complete graph and hypergraphs -- in \cite{balogh} and \cite{west}. Studying size Ramsey numbers for paths has much longer history. Let us mention    
the  breakthrough result by Beck \cite{beck2} that $\sr(P_n,P_n)$  is linear. This result implies that also $\rr(P_n,P_n)$ is linear,
as well as $\rr(P_k,P_n)$ for a fixed $k$, since $\rr(P_k,P_n)\le  \sr(P_k,P_n)$ and we have also an easy general bound $\rr(G,H)\ge |E(G)|+|E(H)|-1$. 
However, it seems not easy to find a multiplicative constant $c$ (if exists) such that $\rr(P_k,P_n)=cn+o(n)$.
In general, we have 
$\rr(P_k,P_n)\le 2n+2k-7$ for every $k,n\ge 2$, proved by Grytczuk, Kierstead and Pra\l at \cite{gryt}, while the best lower bound for $n\ge k\ge 5$ is $\rr(P_k,P_n)\ge 3n/2+k/2-7/2$ by Cyman, Dzido, Lapinskas and Lo \cite{lo}. Thus we know that $3n/2+o(n)\le \rr(P_k,P_n)\le 2n+o(n)$ for $k$ fixed and $n\to \infty$. The authors of   \cite{lo} posed the following conjecture.

\begin{conj}[\cite{lo}]
For every $k\ge 5$
$$\frac{\rr(P_k,P_n)}{n}\to \frac32\quad\text{if } n\to\infty.$$
\end{conj}

We make a step towards this conjecture and prove that $\rr(P_k,P_n)\le 5n/3+o(n)$. 

\begin{theorem}\label{mainpk}
Let $n,k\in\N$ and $k\ge 5$. Then
$$\rr(P_k,P_n)\le\frac53 n+12k.$$  
\end{theorem}
We made no effort to optimise the constant 12 in this theorem. Section \ref{secpk} contains the proof. 

There are very few exact results for $\rr(P_k, P_n)$. It is known that $\rr(P_3,P_n)=\lceil 5(n-1)/4\rceil$ for $n\ge 3$  (\cite{lo})
and there are computer calculated numbers $\rr(P_k,P_n)$ for all $k,n\le 9$ by Pra\l{}at \cite{pral3}. 
It is proved in \cite{lo} that $7n/5-1\le \rr(P_4,P_n)\le 7n/5+9$ for $n\ge 4$ and the authors conjectured that the lower bound is tight.

\begin{conj}[\cite{lo}]
For every $n\ge 4$
$$\rr(P_4,P_n)= \Big\lc \frac75 n\Big\rc -1.$$
\end{conj}

In view of the above mentioned computer calculation, the conjecture is true for  $n\le 9$. We prove that it is true also for $n\ge 10$, by improving the upper bound in \cite{lo}.  

\begin{theorem}\label{mainp4}
For every $n\ge 10$
$$\rr(P_4,P_n)\le \Big\lc \frac75 n\Big\rc -1.$$ 
\end{theorem}

We prove this theorem in Section \ref{secp4}. Our argument is inductive, technically quite complicated, and it is quite different from the argument in \cite{lo}.  
Lately Theorem \ref{mainp4} has been proved independently by Yanbo Zhang and Yixin Zhang \cite{zhang}.

\section{Preliminaries}

For a graph $H=(V,E)$, we put $v(H)=|V|$ and $e(H)=|E|$. 
By the sum of graphs $G$ and $G'$ we mean the union $G\cup G'=(V(G)\cup V(G'), E(G)\cup E(G'))$.

We say that a graph $H$ is \emph{coloured} if every its edge is blue or red. A graph is red (or blue) if all its edges are red (blue). 
We assume that also $\emptyset$ is a coloured graph, which is somewhat non-standard. 
Thus it may happen that  a coloured graph or a subgraph of a coloured graph has 0 vertices.

Let $H$ be a coloured graph. 
By \emph{a component} of $H$ we mean a maximal (in sense of inclusion) connected 
coloured subgraph contained in $H$. If $A,B\subs V(H)$, then $E_H(A,B)$ denotes the set of all
coloured edges of $H$ with one end in $A$ and the other end in $B$.

Given $n\ge k\ge 2$, and a coloured graph $H$ (it may be empty), consider the following auxiliary game $\RR_H(P_k,P_n)$. 
The board of the game is $K_\N$, with exactly $e(H)$ edges coloured and these edges induce a copy of $H$ in $K_\N$. 
The rules of selecting and colouring edges by Builder and Painter are the same as in the standard game $\R(P_k,P_n)$, 
however, Painter is not allowed to colour an edge red if that would create a red $P_k$. 
Builder wins $\RR_H(P_k,P_n)$ at the moment a blue copy of $P_n$ appears on the board.
It is not hard to observe that if Builder has a strategy such that $\RR_\emptyset(P_k,P_n)$ ends within at most $t$ rounds,
then Builder in $\R(P_k,P_n)$ can apply such a strategy as well and finish the game within at most $t$ rounds.

After every round of  $\RR_H(P_k,P_n)$,  
the coloured graph induced by all edges coloured in the game so far (including the edges of $H$) is called \emph{the host graph}. We say a vertex of $K_\N$ is \emph{free} in a round of the game if it is not a vertex of the host graph.

We say that Builder \emph{forces} an edge $uw$ blue if he selects $uw$ 
and Painter has to colour it blue according to the rules of the game, i.e. one of $u,w$, say $u$, is an end of a red path 
on $k-1$ vertices, all distinct from $w$.

\section{Proof of Theorem \ref{mainpk}}\label{secpk}

Let $n,k\in \N$ and $k\ge 5$. It is enough to show a strategy for Builder in $\RR_\emptyset(P_k,P_n)$ 
such that the game ends after at most $5n/3+12 k$ rounds,

In order to simplify the description of Builder's strategy, we assume that Builder can select an edge already coloured. 
In such a round Painter ``colours'' it with the same colour the edge already had. Clearly allowing such moves
cannot help Builder and may only increase the length of the game. 
 
We divide the game $\RR_\emptyset(P_k,P_{n})$ into three stages. Roughly speaking, 
in the first stage Builder creates many blue paths on 3 vertices. 
In the second stage he connects them into at most $k-1$ longer paths, while in the last stage 
he connects a small number of  blue paths into a blue path $P_n$. In order to simplify the description of all three stages, we present a few lemmata.   

\begin{lemma}\label{shortpaths}
Let $n,k\in\N$ and $T=\lc n/3\rc +k$.
Builder has a strategy in $\RR_\emptyset(P_k,P_{n})$ such that
after at most $3T+2(k-1)$ rounds 
the host graph contains $T$ blue, vertex-disjoint paths of length 2.
\end{lemma}

\begin{proof}

We will describe the strategy of Builder based on a definition of active and inactive edges. After every round we will call every coloured edge either active or inactive. An active edge may become inactive
after a few rounds but an inactive edge stays inactive forever.  Given a round of the game, the coloured graph 
induced by all active edges is called the active graph, while by the inactive graph we mean the coloured graph 
induced by all inactive edges. Here is the inductive definition of active and inactive edges.

Before the first round there are no active nor inactive edges and the active and inactive graphs are empty sets. 
Let $t\ge 0$ and suppose that after $t$ rounds we have the active graph $A$ and the inactive 
graph $I$. In the next round Builder chooses an end $x$ of a longest red path $P$ in the active graph
(if $A=\emptyset$, then $x$ is any free vertex) and selects $xy$, where $y$ is any free vertex of the board.   

Suppose Painter colours $xy$ red. It means that $P$ had less than $k-2$ vertices, otherwise a red $P_k$ would have appeared. Then the red edge $xy$ becomes active. After the $(t+1)$-th round the set of active edges is $E(A)\cup\{xy\}$ and 
$E(I)$ is the set of inactive edges.

Suppose Painter colours $xy$ blue. Then we have two possibilities. 
First assume there is a blue edge $e\neq xy$ incident to $x$.
Then the edges $xy$ and $e$ become inactive. If $P$ has a positive length, then the red edge $x'x$, with some $x'\in V(P)$,
also becomes inactive. Let $S=\{xy, e, x'x\}$ if $P$ has a positive length; otherwise let $S=\{xy, e\}$. 
Hence after the $(t+1)$-th round the set of active edges is $E(A)\setm S$, while 
$E(I)\cup S$ is the set of inactive edges. 

Assume now that there is no blue edge $e\neq xy$ incident to $x$. Then the edge $xy$ becomes active. 
After the $(t+1)$-th round the set of active edges is $E(A)\cup\{xy\}$ and 
$E(I)$ is the set of inactive edges.

Builder continues selecting edges in the above way until the graph induced by all blue inactive edges has at least $n+3k$ vertices, then he stops. Let us verify that such a strategy satisfies the assertion of the lemma. 

A routine inductive argument implies that after every round of the game the following holds.

\begin{proposition}\label{stage1}
After every round (until Builder stops the game) the active graph $A$ and the inactive graph $I$ satisfies the following conditions.
\begin{enumerate}[(i)]
\item\label{aidisj}
$E(A)\cap E(I)=\emptyset$.
\item\label{aidisjb}
No inactive blue edge has an endpoint in $V(A)$.
\item\label{aia}
$A$ is a sum of a red path $P$ on less than $k$ vertices and a blue matching $M$ such that every edge of $M$ has exactly one end in $V(P)$.
\item\label{aii}
$I$ is the sum of vertex-disjoint blue paths of length 2 and at most $m$ red edges, where  $m$ is the number of these paths.
\end{enumerate}     
  \end{proposition}

It follows from the above proposition that the game ends at the moment that the inactive graph $I$ contains exactly 
$\lc (n+3k)/3\rc=\lc n/3\rc +k=T$ blue, vertex-disjoint  paths of length 2 and at most $T$ red edges. At this moment 
the active graph $A$ has less than $2(k-1)$ edges, in view of part (\ref{aia}) of the proposition. Thus Builder obtains the required $T$ blue paths within less than $3T+2(k-1)$ rounds. 
\end{proof}

The following definition will be useful in two next lemmata. Given $s,d,m\in\N$ with $d\le s\le m$, we say that a coloured graph $F$ is \emph{an essential $(s,d,m)$-graph}, 
 if it contains $s$ vertex-disjoint paths $G_1,G_2,\ldots,G_s$ satisfying the following conditions.
\begin{enumerate}[(i)]
\item
The paths $G_1,G_2,\ldots,G_s\neq\emptyset$ are blue.    
\item
$F$ contains a red path $P=u_1u_2\ldots u_d$ on $d$ vertices  such that $u_i$ is 
an end of the path $G_i$, for $i=1,2,\ldots d$.
\item
$\sum_{i=1}^s v(G_i)=m$.
\end{enumerate}   
Then we say that $F$ has \emph{the essential red path} $P$ and $s$ \emph{essential blue paths} $G_1,\ldots,G_s$.

\begin{lemma}\label{gluepk1}
Suppose that $1\le d<k\le s\le m$ and after some rounds of the game there is an essential $(s,d,m)$-graph present on the board, 
its essential blue paths are $G_1,G_2,\ldots,G_s$ and $P=u_1u_2\ldots u_d$ is its essential red path $P$. 
Assume that in the next round Builder selects the edge $u_du_{d+1}$, where $u_{d+1}$ is an end of the path $G_{d+1}$.
Then after every response of Painter there is an essential $(s',d',m)$-graph on the board such that 
either $d'=d+1$ and $s'=s$, or $s'= s-1$ and $d'\in\{d,d-1\}$. 
\end{lemma}

\begin{proof}
Let us consider two possible situations after colouring $u_d u_{d+1}$ by Painter.
If $u_d u_{d+1}$ is red, then $P'=P\cup\{u_d u_{d+1}\}$ is a red path on $d'=d+1$ vertices.
It is not hard to verify that the sum of $s$ blue paths $G_i$ and the red path $P'$ is 
an essential $(s',d',m)$-graph that satisfies the required conditions.

Suppose that $u_d u_{d+1}$ is blue. Then we define
blue paths $G'_i=G_i$ for $i=1,2,\ldots,d-1$,
$G'_d=G_d\cup G_{d+1} \cup\{u_d u_{d+1}\}$ and $G'_i=G_{i+1}$ for $i=d+1,d+2.\ldots,s-1$.
We also define ends $u'_i$ of the paths $G'_i$ such that $u'_i=u_i$ for $i=1,2,\ldots,d-1$ and
$u'_i=u_{i+1}$ for $i=d+1,d+2.\ldots,s-1$, while $u'_d$ is one of the ends of $G'_d$.
As for the red path, we put $P'=u_1u_2\ldots u_{d-1}$ if $d\ge 2$; otherwise $P'=u'_1$. 
Thus the red path $P'$ has $d=1$ or $d-1$ vertices, we have $s-1$ vertex-disjoint blue paths $G'_i$ and 
$\sum_{i=1}^s v(G'_i)=\sum_{i=1}^s v(G_i)$. Thus the sum of the red path $P'$ and the blue paths $G'_1,\ldots,G'_{s-1}$
satisfies the required conditions.
\end{proof}

\begin{lemma}\label{lesspaths}
Let $n,k,T\in\N$ and $G$ be a coloured graph containing $T$ blue, vertex-disjoint paths of length 2.
Builder has a strategy in $\RR_G(P_k,P_{n})$ such that
after at most $2T-k$ rounds 
the host graph contains less than $k$ blue, vertex-disjoint paths on $3T$ vertices in total. 
\end{lemma}

\begin{proof}
Let $B\subs G$ be the sum of $T$ blue, vertex-disjoint paths of length 2. Then $B$ is a $(T,1,3T)$-essential graph, with the essential red path of length 0, 
consisting of an end of the first path of $B$.  
Without loss of generality we assume that $G=B$. 

If $T<k$, the Builder in $\RR_B(P_k,P_{n})$ achieves his goal without making any move. Otherwise, Builder selects the edges according to Lemma \ref{gluepk1}, as long as the number of essential blue paths
is at least $k$.  In view of this lemma, we can have rounds of two kinds: when 
the length of the essential red path increases (and the number of essential blue paths does not grow) or when the number of essential blue paths decreases (while the length of the essential red path changes by at most 1). 
A round of the first kind will be called a red-round while a round of the second kind -- a blue-round.    

We will prove that after $2T-k$ rounds or sooner the assertion of the lemma holds. 
Assume for a contradiction that at the end of the $(2T-k)$-th round the number of essential 
blue paths is at least $k$. It means that there where at most $T-k$ rounds such that the number of essential blue paths decreased. Hence the number of rounds such that the length of the essential red path increased was at least $T$. 
Thus the length of the essential red path at the end of the $(2T-k)$-th round is at least $T-(T-k)=k$, which
contradicts the rules of $\RR_G(P_k,P_n)$.

Thus after at most $2T-k$ rounds the host graph contains a sum  $B'$ of less than $k$ vertex-disjoint, blue paths such that $v(B')=v(B)=3T$, what follows from Lemma \ref{gluepk1}. 
\end{proof}

Next two lemmata, useful in the analysis of the third stage, need a definition of a fence. 

Given $s,d,m\in\N$ with $s\le m$, a coloured graph $F$ is called \emph{an $(s,d,m)$-fence},
if it contains $s$ vertex-disjoint paths $G_1,G_2,\ldots,G_s$ satisfying the following conditions.
\begin{enumerate}[(i)]
\item
The path $G_1=w_1\ldots w_d w_{d+1}\ldots w_j$ has at least one vertex, the path  
$P=w_1\ldots w_d$ is red, while the path $G_0=w_d w_{d+1}\ldots w_j$ is blue.
\item
Paths $G_2,\ldots,G_s\neq\emptyset$ are blue.    
\item
$v(G_0)+\sum_{i=2}^s v(G_i)=m$.
\end{enumerate}  
Then we say that $F$ has  \emph{the red picket} $P$ and $s$ \emph{blue pickets} $G_0,G_2,\ldots,G_s$.

\begin{lemma}\label{gluepk2}
Suppose that $m,s\ge 2$, $d\ge 1$, and there is 
an $(s,d,m)$-fence on the board.  
Then Builder can play in such a way that after two rounds the host graph contains 
an $(s',d',m')$-fence such that $m'\ge m-1$ and
either $d'\ge d+1$ and $s'\le s$, or $s'\le s-1$ and $d'\ge d-1$. 
\end{lemma}

\begin{proof}
Suppose there is an $(s,d,m)$-fence on the board, with the red picket 
$P=w_1\ldots w_d$ and blue pickets $G_0,G_2,\ldots,G_s$. Assume that
$w_d$ is an end of $G_0$, while $u_2,u'_2$ are the ends of $G_2$. 
Builder selects the edge $w_du_2$. 

First assume that Painter colours $w_du_2$ red.
Then in the next round Builder selects the same edge $w_du_2$.
After this pair of rounds we define a red path $P'=w_1\ldots w_d u_2$ and the blue paths 
$G'_0=G_2$, $G'_2=G_0\setm\{w_d\}$, $G'_i=G_i$ for $i=3,4,\ldots,s$. 
The paths $G'_1=G'_0\cup P', G'_2,G'_3,\ldots,G'_s$ are vertex-disjoint,
the red path $P'$ has $d+1$ vertices, $v(G'_0)+\sum_{i=2}^s v(G'_i)=m-1$, and 
among $G'_0,G'_2,G'_3,\ldots,G'_s$ only $G'_2$ may be the empty set.
So $\bigcup_{i\le s:\,G'_i\neq\emptyset} G'_i$ is an  $(s',d+1,m-1)$-fence with $s'\in\{s,s-1\}$.

Now assume that Painter colours $w_d u_2$ blue.
Then in the next round Builder selects the edge $w_{d-1}u'_2$, provided $d\ge 2$
(if $d=1$, then Builder selects the same edge  $w_du_2$). Painter colours it red or blue.
After thise pair of rounds we define a red path $P'$ in the following way:
If $d=1$, then $P'=u'_2$; otherwise if $w_{d-1}u'_2$ is red, then $P'=w_1\ldots w_{d-1}u'_2$,
while if $w_{d-1}u'_2$ is blue, then $P'=w_1\ldots w_{d-1}$. We also define blue paths
$G'_i=G_i$ for $i=3,4,\ldots,s$ and a path $G'_0$: if $w_{d-1}u'_2$ is red or $d=1$, then $G'_0=G_0\cup G_2\cup\{w_{d}u_2\}$;
if $w_{d-1}u'_2$ is blue, then $G'_0=G_0\cup G_2\cup\{w_{d}u_2, w_{d-1}u'_2\}$. 
Notice that the paths $G'_1=G'_0\cup P', G'_3,G'_4,\ldots,G'_s$ are vertex-disjoint,
the red path $P'$ has $d$ or $d-1$ vertices, $v(G'_0)+\sum_{i=3}^s v(G'_i)\in\{m,m+1\}$, 
and every path $G'_0,G'_2,G'_3,\ldots,G'_s$ has at least one vertex.
So $G'_1+\bigcup_{i=3}^s G'_i$ is an $(s-1,d',m')$-fence with $d'\le d$ and $m'\ge m$.

In both cases the assertion follows.
\end{proof}

\begin{lemma}\label{onepath}
Let $n,k,t\in\N$ and $G'$ be a coloured graph containing a sum $B$ of $t$ blue, vertex-disjoint paths.
Then Builder has a strategy in $\RR_{G'}(P_k,P_{n})$ such that
after at most $4t+2k$ rounds 
the host graph contains a blue path on at least $v(B)-2t-k$ vertices.
\end{lemma}

\begin{proof}
Without loss of generality we assume that $G'=B$. 
Observe than $B$ consisting of $t$ blue, vertex-disjoint paths is a $(t,1,v(B))$-fence, provided 
we define the red picket of length 0 to  
consist of an end of the first path of $B$.  

Builder in $\RR_{B}(P_k,P_{n})$ selects the edges according to Lemma \ref{gluepk2}, as long as the number of 
blue pickets is at least $2$.
In view of this lemma, we have pairs of rounds of two kinds: when 
the length of the red picket increases (and the number of blue pickets does not grow) or when the number of blue pickets decreases (while the length of the red picket decreases by at most 1). 
A pair of rounds of the first kind will be called a red pair of rounds while a pair of rounds of the second kind -- a blue-pair of rounds.  After every pair of rounds (red or blue) the sum of the number of vertices in blue pickets decreases by
at most one.  

We will prove that after at most $2t+k$ pairs of rounds Builder achieves his goal.  
Assume for a contradiction that at the end of $(2t+k)$-th pair of rounds  
there are at least two blue pickets. 
It means that there where at most $t-2$ blue-pairs of rounds. Hence the number of red-pairs of rounds 
was greater than $t+k$. Thus the length of the red picket is greater than $(t+k)-(t-2)>k$, which
contradicts the rules of $\RR_{G'}(P_k,P_n)$.

Thus after at most $2(2t+k)=4t+2k$ rounds the host graph contains only one blue picket and, in view of Lemma \ref{gluepk2}, the number of its vertices is at least $v(B)-(2t+k)$. 
\end{proof}

We are ready to prove the main theorem.

\begin{proof}[Proof of Theorem \ref{mainpk}]

We present Builder's strategy in the game $\RR_{\emptyset}(P_k,P_n)$ in three stages. 
Let $T=\lc n/3\rc +k$.

\smallskip\noindent
\textbf{Stage 1.} \smallskip

Builder plays according to a strategy whose existence is guaranteed by Lemma \ref{shortpaths} and within at most 
$3T+2(k-1)$ rounds he obtains 
a coloured graph $G$ containing $T$ blue, vertex-disjoint paths of length 2.
The game proceeds to the second stage, equivalent to the game   
$\RR_G(P_k,P_{n})$. 

\smallskip\noindent
\textbf{Stage 2.} \smallskip

In this stage Builder selects the edges according to Lemma \ref{lesspaths} and within at most $2T-k$ rounds of Stage 2
he obtains a coloured graph $G'$ containing a sum $B$ 
of $t<k$ blue, vertex-disjoint paths, such that $v(B)=3T$.  
The game proceeds to the third stage, equivalent to the game   
$\RR_{G'}(P_k,P_{n})$. 

\smallskip\noindent
\textbf{Stage 3.} \smallskip

In the last stage Builder applies a strategy from Lemma \ref{onepath} and after at most $4t+2k$ rounds of Stage 3 the host graph contains a blue path $P$ on at least $v(B)-2t-k$ vertices.
Then Stage 3 ends.

\medskip
Let us analyse the host graph after all three stages of the game. For the blue path $P$ obtained at the end of Stage 3 
we have
$$
v(P)\ge v(B)-2t-k>v(B)-3k=3T-3k\ge n+3k-3k=n.
$$
The number of rounds in all stages is not greater than
$$
(3T+2(k-1))+(2T-k)+(4t+2k)\le 5T+7k-6<5\Big(\frac{n}{3}+k\Big)+7k=\frac53 n+12k.
$$
So a blue path on $n$ vertices was created within at most $5n/3+12k$ rounds and the proof of Theorem \ref{mainpk} is complete.
\end{proof}

\section{Tools for studying the $P_4$ versus $P_n$ game}

Before the proof of Theorem \ref{mainp4} we need a few additional definitions and lemmata. 

If $G$ and $G'$ are coloured graphs and $G'\subs G$, then we denote by $\rn_G(G')$ and $\bn_G(G')$ 
the sets of all red edges and blue edges of $G$ respectively, with at least one end in $V(G')$. 
After every round of a game  $\RR_H(P_4,P_n)$, 
if $G$ is the host graph, for every coloured graph $G'\subs G$ we define  $\rn(G')=\rn_G(G')$ and $\bn(G')=\bn_G(G')$.

Let $c_1,c_2,\ldots,c_k\in\{b,r\}$ be consecutive edge colours of a coloured path $P_{k+1}$. 
Then the coloured path is called a $c_1c_2\ldots c_k$-path. 
Suppose $P\subs K_\N$ is a blue path with ends $x_1,x_2$, and there is a red path of length $l_i\ge 0$ with an end $x_i$, 
for $i=1,2$. Then $P$ is called a blue $(l_1,l_2)$-path and $x_i$ is called its $l_i$-end.   
The coloured path $u_1u_2\ldots u_k u_{k+1}\subs K_\N$ such that $k\ge 2$, 
the edge $u_ku_{k+1}$ is red and $u_1u_2\ldots u_k$
is a blue $(2,0)$-path with a 2-end $u_k$, will be called \emph{an extended $(2,0)$-path}
with \emph{the blue end} $u_1$, \emph{the transition vertex} $u_{k}$ and \emph{the red end} $u_{k+1}$.

Here are two examples of coloured graphs containing a $(1,1)$-path $P$
with $V(P)=\{u_1,u_2,u_3,u_4,u_5\}$ and 1-ends $u_1,u_5$:

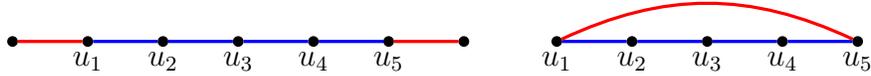
\begin{center}
 \begin{tikzpicture}
\bbbb;
\node [below] at (a1) {$u_1$};
\node [below] at (a2) {$u_2$};
\node [below] at (a3) {$u_3$};
\node [below] at (a4) {$u_4$};
\node [below] at (a5) {$u_5$};
\ver{y}{$(5,0)$};
\ver{x}{$(-1,0)$};
\edr{a5}{y};\edr{a1}{x};
 \end{tikzpicture}\qquad
 \begin{tikzpicture}
\bbbb;
\node [below] at (a1) {$u_1$};
\node [below] at (a2) {$u_2$};
\node [below] at (a3) {$u_3$};
\node [below] at (a4) {$u_4$};
\node [below] at (a5) {$u_5$};
\edar{a5}{a1};
 \end{tikzpicture}
\end{center}

Below there are two examples of coloured graphs containing an extended $(2,0)$-path $P$
on the vertex set $V(P)=\{u_1,u_2,u_3,u_4,u_5,u_6\}$, with the blue end $u_1$, the red end $u_6$ and the transition vertex $u_5$:

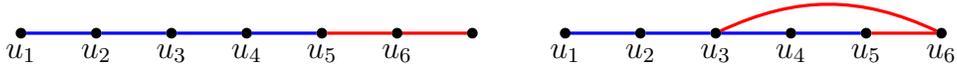
\begin{center}
 \begin{tikzpicture}
\bbbb;
\node [below] at (a1) {$u_1$};
\node [below] at (a2) {$u_2$};
\node [below] at (a3) {$u_3$};
\node [below] at (a4) {$u_4$};
\node [below] at (a5) {$u_5$};
\ver{a6}{$(5,0)$};
\ver{x}{$(6,0)$};
\node [below] at (a6) {$u_6$};
\edr{a5}{a6};\edr{a6}{x};
 \end{tikzpicture}\qquad
 \begin{tikzpicture}
\bbbb;
\node [below] at (a1) {$u_1$};
\node [below] at (a2) {$u_2$};
\node [below] at (a3) {$u_3$};
\node [below] at (a4) {$u_4$};
\node [below] at (a5) {$u_5$};
\ver{a6}{$(5,0)$};
\node [below] at (a6) {$u_6$};
\edr{a5}{a6};\edar{a6}{a3};
 \end{tikzpicture}
\end{center}

The following coloured graph with one red and three blue edges will be called \emph{a limb}.

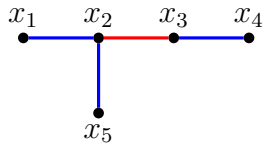
\begin{center}
 \begin{tikzpicture}
\limb;
\node [above] at (l1) {$x_1$};
\node [above] at (l2) {$x_2$};
\node [above] at (l3) {$x_3$};
\node [above] at (l4) {$x_4$};
\node [below] at (l5) {$x_5$};
 \end{tikzpicture}
\end{center}

Suppose that $H$ is a coloured graph that contains vertex-disjoint coloured subgraphs $G_0,G_1,G_2,L$ satisfying the following conditions.

\begin{enumerate}[(A)]
\item\label{limb}
Either $L=\emptyset$, or $L$ is a limb and it is a component of $H$.  
\item\label{path0}
Either $V(G_0)=\emptyset$, or  $G_0$ is a blue path on at least one vertex such that
$G_0\cup\rn_H(G_0)$ is a component of $G$.
Furthermore 
$|\rn_H(G_0)|\le \big\lc\frac25 v(G_0)\big\rc-1$. 
\item\label{path12}
$G_2=\emptyset$ or $G_2$ is an extended $(2,0)$-path on at least three vertices.
Either $G_1=\emptyset$, or $G_1$ is a blue $(1,1)$-path on at least one vertex
such that neither of its 1-ends is adjacent in $G$ to any of: the blue end of $G_2$, the transition vertex of $G_2$, the red end of $G_2$. 
Furthermore $\bn_H(G_1\cup G_2)\setm E(G_1\cup G_2)=\emptyset$
and $|\rn_H(G_1\cup G_2)|\le \big\lc\frac25 v(G_1\cup G_2)\big\rc$. 
\item\label{blue}
The set of all blue edges of $H$ is equal to $\bn_H(G_0\cup G_1\cup G_2\cup L)$. 
\item\label{red}
The set of all red edges of $H$ is equal to $\rn_H(G_0\cup G_1\cup G_2\cup L)$. 
\item\label{div}
$5\mid v(G_0)$ or $5\mid v(G_1\cup G_2)$.
\end{enumerate}
Then we call $H$ \emph{a good graph} with \emph{essential subgraphs} $(G_0,G_1,G_2,L)$. 
We say that $V(L)\cup V(G_0)\cup V(G_1)\cup V(G_2)$ 
is \emph{the set of essential vertices} of $H$ and denote the number of essential vertices of $H$ by $\ess(H)$.
A good graph is \emph{very good} if it satisfies the additional
condition

\begin{enumerate}[(G)]
\item\label{div2}
$5\mid v(G_0)$ and $5\mid v(G_1\cup G_2)$.
\end{enumerate}
 
By \emph{a simple good graph} we mean a good graph such that at least three of its 
essential subgraphs are empty sets. 
Similarly we define \emph{a simple very good graph} $H$, if
$H$ is a very good graph.

Below we present some observations, following immediately from the definition of a good graph, which will be often used, sometimes implicitly, in the paper. 

\begin{proposition}\label{observ}
Suppose that $k\in\N$ and $H$ is a good graph with essential subgraphs $(G_0,G_1,G_2,L)$.
Then the following holds.

\begin{enumerate}[(i)]
\item
The number of red edges of $H$ is not greater than $\big\lc\frac25 \ess(H)\big\rc$ and, if $G_0\neq \emptyset$, 
not greater than $\big\lc\frac25 \ess(H)\big\rc-1$. The number of blue edges of $H$ is not greater than $\ess(H)-1$.
\item\label{goode}
$e(H)\le\big\lc\frac75 \ess(H)\big\rc-1$.
\item\label{del}
$H\setm G_0$ is a good graph, also $H\setm L$  is good.
If $H$ is very good, then also $H\setm G_0$ and $H\setm L$ are very good.
\item\label{verygood}
If $5\mid v(H)$, then $H$ is very good.
\item\label{simple0}
If $H'$ is a sum of a blue path $P$ on $k\ge 1$ vertices and 
at most $\big\lc\frac25 k\big\rc-1$ red edges incident to the path, then $H'$ is a simple good graph, 
with $k$ essential vertices and  essential subgraphs $(P,\emptyset,\emptyset,\emptyset)$.
Furthermore, if $5\mid k$, then $H'$ is simple very good.
\item\label{simple1}
If $H'$ is a sum of a blue $(1,1)$-path $P$ on $k\ge 1$ vertices and 
at most $\big\lc\frac25 k\big\rc$ red edges incident to the path, then $H'$ is a simple good 
graph, with $k$ essential vertices and  essential subgraphs $(\emptyset,P,\emptyset,\emptyset)$.
Furthermore, if $5\mid k$, then $H'$ is simple very good.
\item\label{simple2}
If $H'$ is a sum of an extended $(2,0)$-path $P$ on $k\ge 3$ vertices and 
at most $\big\lc\frac25 k\big\rc$ red edges incident to the path, then $H'$ is a simple good 
graph, with $k$ essential vertices and  essential subgraphs $(\emptyset,\emptyset,P,\emptyset)$.
Furthermore, if $5\mid k$, then $H'$ is simple very good.
\end{enumerate}
\end{proposition}

We omit the uncomplicated justification of Proposition \ref{observ}. The following lemma says, among other things,  that two disjoint very good graphs can be efficiently connected by Builder, if one of them is simple and different from a limb.

\begin{lemma}\label{glue}
Suppose that
$H$ is a very good graph 
and $H'$ is a simple good graph with its essential subgraph
different from a limb.
Then Builder has a strategy in $\RR_{H\cup H'}(P_4,P_{n})$ such that after a finite number of rounds the graph $G$ induced
by all coloured edges of the board is a good graph with $\ess(H)+\ess(H')$ vertices. 

Furthermore, if $H'$ is a simple very good graph, then 
the resulting graph $G$ is very good.
\end{lemma}

\begin{proof}
Let $(G_0,G_1,G_2,L)$ be essential subgraphs of $H$. 
Since $H$ is very good, we have 
$5\mid v(G_0)$ and $5\mid v(G_1\cup G_2)$. 

We split the argument into a few parts, depending on the essential subgraph of $H'$.
If $H'=\emptyset$, the assertion trivially holds so let $H'\neq\emptyset$.

\smallskip\noindent
\textbf{Case 1.} $H'$ has the essential subgraphs $(\emptyset,P,\emptyset, \emptyset)$.\smallskip

Thus $P$ is a blue path $(1,1)$-path on at least 1 vertex and, in view of Condition (\ref{path12}), we have $|\rn_{H'}(P)|\le \lc 2v(P)/5\rc$. 
If $G_1=\emptyset$ then, since $5\mid v(G_0)$,
clearly $H\cup H'$ is a good graph with essential subgraphs $(G_0,P,G_2,L)$. Thus the required
good graph is obtained without any move in the game. 
 
Assume further that $G_1\neq\emptyset$. Let $x,y$ be the 1-ends of $P$
and $u_1,u'_1$ be the 1-ends of $G_1$. 
 Builder in $\RR_{H\cup H'}(P_4,P_{n})$  forces the edge $xu_1$ blue.
Then a blue $(1,1)$-path $G'_1$ with 1-ends $u'_1,y$ appears on the board and it has 
$v(G_1)+v(P)$ vertices. Let $G$ be the host graph at this moment.
Since $5\mid v(G_1\cup G_2)$, we have
\begin{eqnarray*}
|\rn_G(G'_1\cup G_2)|&=&|\rn_H(G_1\cup G_2)|+|\rn_{H'}(P)|\le \frac25 v(G_1\cup G_2)+\Big\lc \frac25 v(P)\Big\rc\\
&=&\Big\lc \frac25 \big(v(G_1)+v(P)+v(G_2)\big)\Big\rc=\Big\lc \frac25 v(G'_1\cup G_2)\Big\rc.
\end{eqnarray*}
Furthermore 
$\rn_G(G_0)=\rn_H(G_0)$, $\bn_G(G_0)=\bn_H(G_0)$, $\rn_G(L)=\rn_H(L)$, $\bn_G(L)=\bn_H(L)$,  
and $5\mid v(G_0)$ by the assumption that $H$ is very good, so in view of the above estimation of $|\rn_G(G'_1\cup G_2)|$, 
we conclude that $G$ is a good graph with essential subgraphs $(G_0,G'_1,G_2,L)$ and $\ess(G)=\ess(H)+\ess(H')$.

\smallskip\noindent
\textbf{Case 2.} $H'$ has the essential subgraphs $(\emptyset,\emptyset, P,\emptyset)$.\smallskip

Thus $P$ is an extended $(2,0)$-path on at least 3 vertices, with $|\rn_{H}(P)|\le \lc 2v(P)/5\rc$. 
If $G_2=\emptyset$, then $H\cup H'$ is the required good graph, 
with essential subgraphs $(G_0,G_2,P,L)$. Assume further that $G_2\neq\emptyset$.

Let $y,x,z$ be the blue end, the red end and the transition vertex of $P$, respectively.
Similarly, let $u'_2,u_2,w_2$ be the blue end, the red end and the transition vertex of $G_2$.

Builder in $\RR_{H\cup H'}(P_4,P_{n})$  forces the edge $zu'_2$ blue.
Then we obtain an extended $(2,0)$-path $G'_2$ on the vertex set $V(G_2\cup (P\setm\{x\})$, with its blue end $y$ and its red end $u_2$,
and the blue $(1,1)$-path on one vertex $x$. 

If $G_1=\emptyset$, we define $G'_1=x$; otherwise let $u_1,u'_1$ be the 1-ends of $G_1$. 
In the latter case, in the next round Builder forces the edge $xu_1$ blue and we define $G'_1=G_1\cup \{xu_1\}$.
In both cases $G'_1$ is a blue $(1,1)$-path on the vertex set $V(G_1)\cup \{x\}$.
Let us estimate the number of red edges incident to vertices of $G'_1\cup G'_2$ at this moment of the game, 
bearing in mind the assumption $5\mid v(G_1\cup G_2)$. We have  
\begin{eqnarray*}
|\rn(G'_1\cup G'_2)|&=&|\rn_{H'}(P)|+|\rn_H(G_1\cup G_2)|
\le \Big\lc \frac25 v(P)\Big\rc+\frac25 v(G_1\cup G_2)\\
&=&\Big\lc\frac25 v(P\cup G_1\cup G_2)\Big\rc
=\Big\lc\frac25 v(G'_1\cup G'_2)\Big\rc.
\end{eqnarray*}
As in the previous case, adding a blue edge $xu_1$ to the host graph does not change the 
set of edges incident to any vertex of $G_0$ and $L$, and we have $5\mid v(G_0)$, 
thus the obtained host graph is a good graph 
with essential subgraphs $(G_0,G'_1,G'_2,L)$, with $\ess(H)+\ess(H')$ essential vertices.

\smallskip\noindent
\textbf{Case 3.} $H'$ has the essential subgraphs $(P,\emptyset,\emptyset, \emptyset)$.\smallskip

Thus $P$ is a blue path on at least 1 vertex and 
$|\rn_{H'}(P)|\le \lc 2v(P)/5\rc -1$. If $G_0=\emptyset$ then, since $5\mid v(G_1\cup G_2)$,
clearly $H\cup H'$ is a good graph with essential subgraphs $(P,G_1,G_2,L)$. Thus the required
good graph is obtained at once. 
 
Assume further that $G_0\neq\emptyset$, $u_0,u'_0$ are the ends of $G_0$ 
and $x,y$ are the ends of $P$. 
Builder starts $\RR_{H\cup H'}(P_4,P_{n})$  with selecting the edge $xu_0$.

If Painter colours $xu_0$ blue, then we obtain a blue path $G'_0$ on $v(G_0)+v(P)$ vertices and
\begin{eqnarray*}
|\rn(G'_0)|&=&|\rn_H(G_0)|+|\rn_{H'}(P)|\le \frac25 v(G_0)-1+\Big\lc \frac25 v(P)\Big\rc-1\\
&=&\Big\lc \frac25 (v(G_0)+v(P))\Big\rc-2<\Big\lc \frac25 v(G'_0)\Big\rc-1.
\end{eqnarray*}
The above estimation and the fact that adding the blue edge $xu_0$ does not affect 
the set of edges incident to $G_1$, $G_2$ or $L$, implies that for 
$5\mid v(G_1\cup G_2)$, the obtained host graph is a good graph 
with essential subgraphs $(G'_0,G_1,G_2,L)$ and $v(G'_0)+v(G_1\cup G_2)+v(L)=\ess(H)+\ess(H')$
essential vertices.

Now assume that Painter colours  $xu_0$ red. Then Builder selects $yu'_0$.
If Painter colours it blue, then a blue $(1,1)$-path $G'_1$ with 1-ends $u_0,x$ appears on the board.
Suppose Painter colours $yu'_0$ red. Then, if $P$ has more than one vertex,  Builder forces the edge $yu_0$ blue,
and thereby obtains a blue $(1,1)$-path $P'$ with 1-ends $u'_0,x$. If $P$ has only one vertex, we get 
an extended $(2,0)$-path $P'$, with its blue end $u'_0$ and its red end $x$, containing the red edge $u_0x$. 
In all cases the obtained path $P'$ has $v(G_0)+v(P)$ vertices and,  
since at most two new red edges (incident to vertices of $P'$) were selected in the game, we have 
\begin{eqnarray*}
|\rn(P')|&\le& |\rn_H(G_0)|+|\rn_{H'}(P)|+2\le \frac25 v(G_0)-1+\Big\lc \frac25 v(P)\Big\rc-1+2\\
&=&\Big\lc \frac25 (v(G_0)+v(P))\Big\rc=\Big\lc \frac25 v(P')\Big\rc.
\end{eqnarray*}
Let $H''$ be the sum of $P'$ and all (at most two) red edges selected in the game. 
Then the above estimation and Proposition \ref{observ}(\ref{simple1}),(\ref{simple2}) imply that $H''$ is a simple good graph. 
Furthermore, in view of Proposition \ref{observ}(\ref{del}),   
the graph $H\setm G_0$ is very good. 
Let $F$ be the host graph at this moment of the game. 
Observe that $F$ is the sum of $(H\setm G_0)$ and  $H''$ and  all assumptions of
Lemma \ref{glue} are fulfilled by the graphs $H\setm G_0$ and $H''$.
We have already proved Lemma \ref{glue}
if the simple good graph contains an extended $(2,0)$-path (Case 2) or a blue $(1,1)$-path (Case 1)
so we argue that after a few rounds of  $\RR_{F}(P_4,P_n)$ the obtained host graph is a good graph such that its number of essential vertices  is
$$\ess(H\setm G_0)+\ess(H'')=v(G_1\cup G_2)+v(L)+v(P')=v(G_1\cup G_2)+v(L)+v(G_0)+v(P)=\ess(H)+\ess(H').$$
Hence the first part of the assertion follows. 

The second part of the assertion is a consequence of Proposition \ref{observ}(\ref{verygood}) and the fact that if both graphs $H$ and $H'$ are very good, then
$5\mid \ess(H)+\ess(H')$.  
\end{proof}

The next lemma is crucial to the inductive argument in the proof of Theorem \ref{mainp4}.

\begin{lemma}\label{good}
Suppose that $n\ge 10$, $5\le k\le n$, $5\mid k$, and $H$ is a coloured graph.
Assume that $H$ is good and has $k-5$ essential vertices.
Then Builder has a strategy in $\RR_H(P_4,P_{n})$ such that 
after a finite number of rounds the host graph is a very good graph with $k$ essential vertices. 
\end{lemma}

\begin{proof}
Suppose that $H$ satisfies the assumption of the lemma and has essential subgraphs $(G_0,G_1,G_2,L)$.
Note that $H$ is very good, in view of \ref{observ}(\ref{verygood}).

Suppose that $K_\N$ contains the coloured graph $H$ and $K'_\N\subs K_\N$ is the complete graph, vertex-disjoint from $H$.
Builder starts $\RR_H(P_4,P_n)$ by selecting two adjacent edges of $K'_\N$, say $ab$ and $bc$. We consider all Painter's responses.

\smallskip\noindent
\textbf{Case 1.} $ab$ and $bc$ are red.\smallskip

Builder selects edges $ax$ and $cy$ for any free, distinct $x,y\in V(K'_\N)$. 
According to the rules of the game $\RR_H(P_4,P_n)$, Painter has to colour them blue
so the $brrb$-path is obtained:

\begin{center}
 \begin{tikzpicture}
\brrb;
\node [above] at (a1) {$x$};
\node [above] at (a2) {$a$};
\node [above] at (a3) {$b$};
\node [above] at (a4) {$c$};
\node [above] at (a5) {$y$};
 \end{tikzpicture}
\end{center}
Then Builder forces $xc$ blue and the following coloured graph appears:

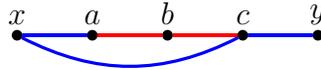
\begin{center}
 \begin{tikzpicture}
\brrb;
\node [above] at (a1) {$x$};
\node [above] at (a2) {$a$};
\node [above] at (a3) {$b$};
\node [above] at (a4) {$c$};
\node [above] at (a5) {$y$};
\edab{a1}{a4};
 \end{tikzpicture}
\end{center}
Thus it contains an extended $(2,0)$-path  $ycxab$.
In view of Proposition \ref{observ}(\ref{simple2}), the coloured graph $H'$ induced by the five edges selected in the game
 is a simple very good graph with five essential vertices. 
Since the host graph at this moment is the sum of vertex-disjoint very good graphs $H$ and $H'$, we can apply Lemma \ref{glue} and  
we infer that Builder can continue the game so that
after some rounds the host graph is a very good graph, with $\ess(H)+\ess(H')=k$ essential vertices.

\smallskip\noindent 
\textbf{Case 2.} $ab$ and $bc$ are blue.\smallskip

In the next two rounds Builder selects edges $cx$ and $xy$ for any free, distinct $x,y\in V(K'_\N)$.
Let us consider four cases, depending on Painter's response. 

If Painter colours them blue, we obtain a blue path $H'=abcxy$ such that no red edge is incident to it. 
Based on Proposition \ref{observ}(\ref{simple0}), it is a simple very good graph on five vertices and again 
we apply Lemma \ref{glue} to $H$ and $H'$. Thereby  
we infer that after some rounds Builder obtains a very good host graph with $\ess(H)+\ess(H')=k$ essential vertices.

If Painter colours $cx$ and $xy$ red, then Builder forces the edge $ay$ blue and
we obtain an extended $(2,0)$-path $cbayx$.

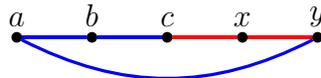
\begin{center}
 \begin{tikzpicture}
\bbrr;
\node [above] at (a1) {$a$};
\node [above] at (a2) {$b$};
\node [above] at (a3) {$c$};
\node [above] at (a4) {$x$};
\node [above] at (a5) {$y$};
\edab{a1}{a5};
 \end{tikzpicture}
\end{center}

The coloured graph $H'$ induced by the five edges selected in the game
 is a simple very good graph with five essential vertices, we apply Lemma \ref{glue} again, and conclude that
after some rounds Builder obtains a very good graph with $\ess(H)+\ess(H')=k'$ essential vertices.
 
Suppose Painter colours $cx$ red and $xy$ blue. So we obtain a $bbrb$-path:

\begin{center}
 \begin{tikzpicture}
\bbrb;
\node [above] at (a1) {$a$};
\node [above] at (a2) {$b$};
\node [above] at (a3) {$c$};
\node [above] at (a4) {$x$};
\node [above] at (a5) {$y$};
 \end{tikzpicture}
\end{center}
Then Builder selects $ay$.
If Painter colours it blue, we receive a blue $(1,1)$-path on five vertices, with 1-ends $x,c$.
If Painter colours $ay$ red, Builder forces $cy$ blue and we receive the following
coloured graph.

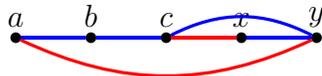
\begin{center}
 \begin{tikzpicture}
\bbrb;
\node [above] at (a1) {$a$};
\node [above] at (a2) {$b$};
\node [above] at (a3) {$c$};
\node [above] at (a4) {$x$};
\node [above] at (a5) {$y$};
\edab{a5}{a3};\edar{a1}{a5};
 \end{tikzpicture}
\end{center}
It is a blue $(1,1)$-path with 1-ends $x,a$.

Thus, regardless of whether Painter coloured $ay$ red or blue, Builder obtains 
 a graph $H'$ which is a sum of a blue $(1,1)$-path on five vertices and at most two red edges
incident to the path.  As before, we apply Proposition \ref{observ}(\ref{simple1}), then Lemma \ref{glue} 
to $H'$ and $H$, and conclude the assertion.

Suppose Painter colours $cx$ blue and $xy$ red. Then Builder selects $ay$.
If Painter colours it blue, we receive a blue $(1,1)$-path on five vertices, with 1-ends $x,y$.
If Painter colours $ay$ red, we receive an extended $(2,0)$-path $abcxy$. 
The analysis is analogous as before and we conclude the assertion in a similar way.

\smallskip\noindent 
\textbf{Case 3.} $ab$ and $bc$ have different colours, i.e.~ $ab$ is, say, blue and $bc$ is red.\smallskip

Suppose first that $L=\emptyset$. 
Builder selects $cx$ for any free $x\in V(K'_\N)$. 

If Painter colours it red, then Builder forces an edge $xy$ blue for any free $y\in V(K'_\N)$, 
then he forces $ax$ blue. The following coloured graph is obtained:

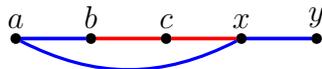
\begin{center}
 \begin{tikzpicture}
\brrb;
\node [above] at (a1) {$a$};
\node [above] at (a2) {$b$};
\node [above] at (a3) {$c$};
\node [above] at (a4) {$x$};
\node [above] at (a5) {$y$};
\edab{a1}{a4};
 \end{tikzpicture}
\end{center}
Such a graph with an extended $(2,0)$-path has been already analysed in Case 1. 

Suppose  Painter colours $cx$ blue. Then Builder selects $by$ for any free $y\in V(K'_\N)$.
If Painter colours it blue, the four edges selected in the game so far form a limb $L'$, vertex-disjoint with $H$.
Hence, since $H$ is very good and $L=\emptyset$, also $H\cup L'$ is a very good graph. It has $\ess(H)+5=k$ 
essential vertices and the assertion follows. 
Assume that Painter  colours $by$ red. We receive the following coloured graph at $K'_\N$:

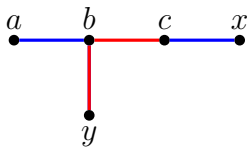
\begin{center}
 \begin{tikzpicture}
\limb;
\node [above] at (l1) {$a$};
\node [above] at (l2) {$b$};
\node [above] at (l3) {$c$};
\node [above] at (l4) {$x$};
\node [below] at (l5) {$y$};
\edr{l2}{l5};
 \end{tikzpicture}
\end{center}
Then Builder forces edges $ac$ and $xy$ blue and 
we obtain a graph $H'$ which is a sum of a blue $(2,1)$-path (which is also a $(1,1)$-path) 
$yxcab$ and a red path of length 2, as we illustrate in the following picture.

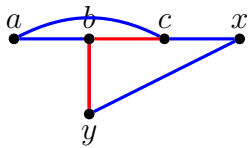
\begin{center}
 \begin{tikzpicture}
\limb;
\node [above] at (l1) {$a$};
\node [above] at (l2) {$b$};
\node [above] at (l3) {$c$};
\node [above] at (l4) {$x$};
\node [below] at (l5) {$y$};
\edr{l2}{l5};
\edab{l3}{l1};\edb{l4}{l5};
 \end{tikzpicture}
\end{center}
Such a graph is simple very good so based on Lemma \ref{glue} applied to $H'$ and $H$, after some rounds 
Builder obtains a very good graph with $\ess(H)+\ess(H')=k$ essential vertices.

Now suppose that $L\neq\emptyset$. 
Then Builder selects an edge $xy$ for any free $x,y\in V(K'_\N)$. 

If Painter colours $xy$ blue, then Builder selects $cx$. No matter how Painter colours it, 
the four edges  selected in the game form either a $brrb$-path considered in
Case 1 or a $bbrb$-path, which we considered in Case 2.
Therefore further we assume that Painter colours $xy$ red. Then Builder selects $ax$. 

If Painter colours $ax$ blue, then Builder forces $by$ and $cy$ blue. The following graph is obtained:

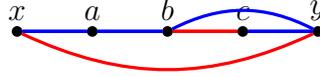
\begin{center}
 \begin{tikzpicture}
\bbrb;
\node [above] at (a1) {$x$};
\node [above] at (a2) {$a$};
\node [above] at (a3) {$b$};
\node [above] at (a4) {$c$};
\node [above] at (a5) {$y$};
\edab{a5}{a3};\edar{a1}{a5};
 \end{tikzpicture}
\end{center}
Such a stage of the game has been already analysed in Case 2.

Further we assume that $ax$ is red. Consider two components of the coloured graph induced 
by all coloured edges present on the board $K_\N$: the $rrbr$-path created in the four rounds of the game and
the limb $L$, with its vertices denoted as in the following picture. 

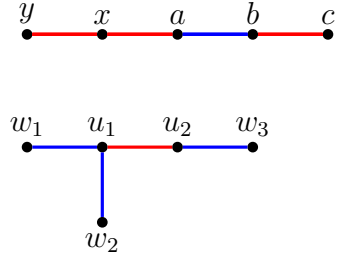
\begin{center}
 \begin{tikzpicture}
\rrbr;
\node [above] at (a1) {$y$};
\node [above] at (a2) {$x$};
\node [above] at (a3) {$a$};
\node [above] at (a4) {$b$};
\node [above] at (a5) {$c$};
\limb;
\node [above] at (l1) {$w_1$};
\node [above] at (l2) {$u_1$};
\node [above] at (l3) {$u_2$};
\node [above] at (l4) {$w_3$};
\node [below] at (l5) {$w_2$};
 \end{tikzpicture}
\end{center} 
In the next five rounds Builder forces blue edges: $yw_1$, $yw_3$, $cu_2$, $cx$, $aw_2$. 
After we obtain the following  coloured graph $H'$:

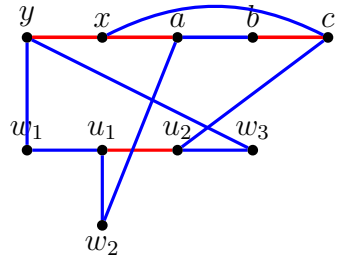
\begin{center}
 \begin{tikzpicture}
\rrbr;
\node [above] at (a1) {$y$};
\node [above] at (a2) {$x$};
\node [above] at (a3) {$a$};
\node [above] at (a4) {$b$};
\node [above] at (a5) {$c$};
\limb;
\node [above] at (l1) {$w_1$};
\node [above] at (l2) {$u_1$};
\node [above] at (l3) {$u_2$};
\node [above] at (l4) {$w_3$};
\node [below] at (l5) {$w_2$};
\edb{a1}{l1};\edb{a1}{l4};\edb{a3}{l5};\edb{a5}{l3};\edab{a5}{a2};
 \end{tikzpicture}
\end{center}
Thus $H'$ is a sum of a blue $(1,1)$-path on 10 vertices and four red edges. 
We know by Proposition \ref{observ}(\ref{simple1}) that such a graph is simple very good, and that, by 
Proposition \ref{observ}(\ref{del}), the graph  $H\setm L$ is very good.
The very good graphs $H'$ and  $H\setm L$ are vertex-disjoint 
and they satisfy the assumption of Lemma \ref{glue}.  Similarly to the previous analysis,
based on Lemma \ref{glue} Builder obtains a very good graph with $\ess(H\setm L)+\ess(H')=k$ essential vertices.
\end{proof}

An immediate consequence of Lemma \ref{good} and the inductive argument 
is the following corollary.
 
\begin{cor}\label{good2}
Suppose that $n\ge 10$. 
Then Builder has a strategy in $\RR_\emptyset(P_4,P_{n})$ such that 
after a finite number of rounds the host graph is a very good graph with $5\lf n/5\rf$ essential vertices. 
\end{cor}

So far Lemma \ref{glue} was applied to very good graphs only.
In the next section we will use it also in case of a simple good graph $H'$ which is not very good, i.e.~with the number of vertices not divisible by 5.

\section{Proof of Theorem \ref{mainp4}}\label{secp4}

Let $n=m+r$ with $5\mid m$, $m\ge 10$  and $0\le r\le 4$. 
 It is enough to show a strategy for Builder in $\RR_\emptyset(P_4,P_n)$ 
such that the game ends after at most $\lc 7n/5\rc-1$ rounds.
We divide the game $\RR_\emptyset(P_4,P_{n})$ into three stages. Roughly speaking, 
in the first stage Builder creates a big very good graph
and in the second stage he increases it to a good graph with $n$ essential vertices. 
In the last stage he connects parts of the good graph into a blue path $P_n$.  

\smallskip\noindent
\textbf{Stage 1.} \smallskip

In the first stage, based on  Corollary \ref{good2}, Builder uses a strategy which guarantees
that after some round the host graph is a very good graph $H$ with $m$ essential vertices. 
Then the first stage ends. 
Assume that the very good graph $H$ has essential subgraphs $(G_0,G_1,G_2,L)$.
The game proceeds to the next stage, equivalent to the game $\RR_H(P_4,P_{n})$.

\smallskip\noindent
\textbf{Stage 2.} \smallskip

Let $K'_\N\subs K_\N$ be a complete graph, vertex-disjoint from $H$.
Builder begins by selecting $r-1$ edges of a path on $r$ vertices in $K'_\N$ (if $r\le 1$, he does nothing). 
After Painter colours them, we obtain a coloured path $P$ on $r$ vertices.
For $r=0$ we put $P=\emptyset$, for $r=1$ we have a trivial path $P$.  
We consider a few cases in order to define a new coloured component $H'$. The only case not listed below is
when $r=2$ and $P$ is a red edge. We call it \emph{the exceptional case}, assume that the second stage ended here 
and consider this case later. 

\smallskip\noindent
\textbf{Case 1.} $0\le r\le 4$ and $P$ is blue.\smallskip

Then we put $H'=P$. There is no red edges incident to the path $H'$ so it is a simple good graph with $\ess(H')=r$.

\smallskip\noindent
\textbf{Case 2.} $r=3$, $P=abc$ and it is an $rr$-path.\smallskip

Then Builder forces blue edges $ax,xc$ with any free $x\in V(K'_\N)$. 
Hence we obtain a coloured graph $H'$ that is the sum of a blue $(1,1)$-path $G_1=axc$ and two red edges.
$H'$ is a simple good graph since 
$|\rn(H')|=2=\lc 2 v(G_1)/5\rc$.  
Clearly $\ess(H')=r$.

\smallskip\noindent
\textbf{Case 3.} $r=3$, $P=abc$ and it is a $br$-path.\smallskip

Then Builder selects $ac$. No matter how Painter colours it, the coloured graph induced by 
the three edges form a simple good graph $H'$ with 3 essential vertices. Indeed, 
if $ac$ is blue, then $H'$ contains a blue $(1,1)$-path $cab$;
if $ac$ is red, then $H'$ contains an extended $(2,0)$-path $abc$.
  
\smallskip\noindent
\textbf{Case 4.} $r=4$, $P=abcd$ and at least one edge of $P$ is red.\smallskip

Then Builder selects $ad$. Thus the four selected edges form a cycle. 
If three edges of the cycle are blue, we get a blue $(1,1)$-path on 4 vertices.
If two edges of the cycle are blue and two red edges are adjacent, we get an extended $(2,0)$-path on 4 vertices.
If the two red edges are disjoint, say $ab$ and $cd$, then Builder forces a blue edge $ac$ and
we obtain a blue $(1,1)$-path on 4 vertices. 
Obviously, the cycle cannot have exactly one edge blue, since Painter never creates a red $P_4$. 
Thus in all cases for $r=4$ the four coloured edges form a simple good graph and we denote it by $H'$. 
We have $\ess(H')=4$.

It follows from the argument above that for every $r\le 4$, in every case apart from the exceptional case 
($r=2$ and $P$ red) the new component $H'$ is a simple good 
graph with $r$ essential vertices. Furthermore, $H$ and $H'$ satisfy the assumption of Lemma \ref{glue}.
Thus within some further rounds  Builder creates a good graph $G$ with $\ess(G)=m+r=n$ and
 with essential subgraphs $(G'_0,G'_1,G'_2,L)$. 

The second stage ends and the game proceeds to the third stage, equivalent to the game $\RR_{G}(P_4,P_n)$.

\smallskip\noindent
\textbf{Stage 3.} \smallskip

Let us recall that we still exclude  the exceptional case  from the analysis.
Builder begins by transforming the limb $L$ (if $L\neq\emptyset$) into a blue $(1,0)$-path.
Suppose the vertices of $L$ are denoted as in the following picture.
 
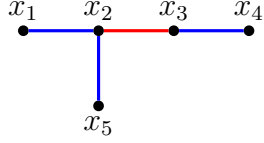
\begin{center}
 \begin{tikzpicture}
\limb;
\node [above] at (l1) {$x_1$};
\node [above] at (l2) {$x_2$};
\node [above] at (l3) {$x_3$};
\node [above] at (l4) {$x_4$};
\node [below] at (l5) {$x_5$};
 \end{tikzpicture}
\end{center}
Builder selects $x_4x_5$. If Painter colours it blue, we obtain a blue $(1,0)$-path 
$P'=x_3x_4x_5x_2x_1$. If Painter colours $x_4x_5$ red, Builder forces the edge $x_5x_3$ blue and we get  a blue $(1,0)$-path 
$P'=x_1x_2x_5x_3x_4$, as shown in the following picture.

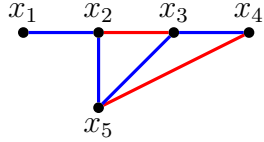
\begin{center}
 \begin{tikzpicture}
\limb;
\node [above] at (l1) {$x_1$};
\node [above] at (l2) {$x_2$};
\node [above] at (l3) {$x_3$};
\node [above] at (l4) {$x_4$};
\node [below] at (l5) {$x_5$};
\edr{l4}{l5};\edb{l3}{l5};
 \end{tikzpicture}
\end{center}

Without loss of generality we can assume that $x_4$ is an 1-end of the path $P'$. 
Clearly $|\rn(P')|\le 2=2 v(P')/5$.   
If $G'_1=\emptyset$, we put $F_1=P'$; otherwise
for an 1-end $u_1$ of $G_1$, Builder forces the edge $x_4u_1$ blue. 
Then we define $F_1=P'\cup G'_1\cup\{x_4u_1\}$. 
 
For $L=\emptyset$ we define $F_1=G'_1$. 
Thus in both cases ($L$ empty or not) 
$F_1$ is a $(1,0)$-path or $F_1=\emptyset$.
Observe also that in both cases at this moment of the game we have   
$v(F_1)=v(L)+v(G'_1)$ and 
$|\rn(L)|=2 v(L)/5$.
 
If $F_1=\emptyset$ or $G'_2=\emptyset$, we define $D=F_1\cup G'_2$. 
Suppose now that $F_1,G'_2\neq\emptyset$. Let $f_1$ be an 1-end of the path $F_1$,
$f'_1$ be its other 1-end, and let $u'_2,w_2,u_2$ denote the blue end, the transition vertex and the red end of $G'_2$, 
respectively. Breaker forces the edges $f_1u_2$ and $w_2f'_1$ blue.
After there is a blue $(1,0)$-path $D$ on $v(F_1)+v(G'_2)$ vertices on the board, 
with an 1-end $u_2$ and the other end $u'_2$. 

In both cases ($F_1,G'_2\neq\emptyset$ or not) 
we have $v(D)=v(F_1)+v(G'_2)=v(L)+v(G'_1)+v(G'_2)$ and every red edge incident to a vertex of $D$
is either a red edge of $G$ incident to $G'_1$ or $G'_2$, or it is incident to $L$.
Thus
\begin{equation}\label{rnd}
|\rn(D)|=|\rn_G(G'_2\cup G'_1)|+|\rn(L)|=|\rn_G(G'_2\cup G'_1)|+\frac25 v(L).
\end{equation}
Furthermore, every blue edge of $L\cup G'_1\cup G'_2$ is a blue edge of $D$.

It remains to connect $D$ and $G'_0$ into a blue path $D'$. 
Let $x$ be an 1-end of $D$ and $x'$ be the other end of it.
If  $G'_0=\emptyset$, then we define $D'=D$.
Clearly then $|\rn(D')|=|\rn(D)|$ and
\begin{eqnarray*}
|\rn(D')|&\le& |\rn_G(G'_2\cup G'_1)|+\frac25 v(L)
\le \Big\lc \frac 25v(G'_1\cup G'_2)\Big\rc+\frac25 v(L)\\
&=&
\Big\lc \frac 25(\ess(G)-v(L))\Big\rc+\frac25 v(L)=
\Big\lc \frac 25 \ess(G)\Big\rc.
\end{eqnarray*}

Suppose that $G'_0\neq\emptyset$. Let $u_0$ and $u'_0$ be the ends of the path $G'_0$. 
Builder selects the edge $xu_0$. 
If Painter colours it blue, then we define $D'=D\cup G'_0\cup\{xu_0\}$. 
Suppose Painter colours $xu_0$ red. Then $u_0$ becomes a 2-end of $G'_0$. 
We also know that there is a red edge, say $e$, incident to $x$ but different from $xu_0$. 
If $e$ is not incident to $x'$, then Builder forces the edge $u_0x'$ blue. 
If $e$ is incident to $x'$, then $x'$ is a 2-end of $D$ now and Builder can force blue
either the edge $x'u'_0$ if $u'_0\neq u_0$, or an edge $x'y$ with a free $y$ if $u'_0= u_0$. 
We denote by $D'$ the resulting blue path on $v(D)+v(G'_0)$ vertices. Since 
at most one red edge incident to its vertex was selected while connecting $D$ and $G'_0$,
in view of \eqref{rnd}  we have
\begin{eqnarray*}
|\rn(D')|&\le& |\rn_G(G'_2\cup G'_1)|+\frac25 v(L)+|\rn_G(G'_0)|+1\\
&\le&
\Big\lc \frac 25v(G'_1\cup G'_2)\Big\rc+\frac25 v(L)+\Big\lc\frac25 v(G'_0)\Big\rc=\Big\lc\frac25 \ess(G)\Big\rc.
\end{eqnarray*}
In the last equality we used the property of a good graph that one of the numbers $v(G'_1\cup G'_2)$, $v(G'_0)$
is divisible by 5 and, clearly, $2 v(L)/5$ is integer. 

In all cases we obtain a blue path $D'$ on $v(L)+v(G'_0)+v(G'_1)+v(G'_2)+v(L)=\ess(G)=n$ vertices,
such that 
$|\rn(D')|\le \lc 2 \ess(G)/5\rc=\lc 2 n/5\rc$. 
The third stage ends.

\medskip
Observe that the set $\rn(D')$ for the path $D'$ obtained at the end of Stage 3 is also the set of all red edges selected 
in the game. Indeed, $G$ was good at the end of Stage 2 so every red edge was incident to 
$G'_0\cup G'_1\cup G'_2\cup L$ then. Furthermore every red edge selected during the third stage 
was incident to $L$ or to the path $D$ with $V(D)=V(L)\cup V(G'_1)\cup V(G'_2)$. 
It follows from the construction of $D'$ that $D\subs D'$. It may have happened that $G'_0\not\subs D'$, but
only if $G'_0$ had one vertex. However, in such a case at the end of Stage 2 the number of red edges
incident to $G'_0$ was at most $\lc 2 v(G'_0)/5\rc-1=0$. Thus the number of red edges selected 
in the game is $|\rn(D')|\le \lc 2 n/5\rc$. 

Notice also that at the end of the third stage there are no blue edges on the board other than the edges of $D'$
so the number of blue edges selected in the game $\RR_\emptyset(P_4,P_n)$ is $n-1$.

Summarising, the number of all coloured edges (and thus the number of 
rounds of the game $\RR_\emptyset(P_4,P_n)$) is not greater than  
$$n-1+\Big\lc\frac25 n\Big\rc=\Big\lc\frac75 n\Big\rc-1.$$
This ends the proof in all cases apart from the exceptional case. 

We proceed to the exceptional case. Let us recall that then we assume that $r=2$, after the second stage 
the position on the board consists of a red path $P=ab$ and the coloured graph $H$ (vertex-disjoint from $P$)
which is very good, with essential subgraphs $(G_0,G_1,G_2,L)$ and $m$ essential vertices. 
The game proceeds to the third stage.

\smallskip\noindent
\textbf{Stage 3 in the exceptional case.} \smallskip

We consider three subcases, depending on the essential subgraphs of $H$.

\smallskip\noindent
\textbf{Subcase 1.} $G_0\neq\emptyset$.\smallskip

Let $u_0,u'_0$ be the ends of $G_0$.
Builder selects edges $u_0a$ and $u'_0b$. Since Painter avoids a red $P_4$, he has to colour blue at least one of these
edges. If both are blue, we obtain a blue $(1,1)$ path $B$ with 1-ends $a,b$, on $v(G_0)+v(P)$ vertices.
If exactly one of the edges is blue, we get an extended $(2,0)$-path on $v(G_0)+2$ vertices.
In both cases, since $5\mid v(G_0)$, the new coloured component $H'=B\cup P\cup\{u_0a,u'_0b\}$ on the board satisfies 
$$|\rn(B)|\le |\rn_H(G_0)|+|\rn(P)|\le \frac25 v(G_0)-1+2=\frac25 v(G_0)+\Big\lc\frac25 v(P)\Big\rc
=\Big\lc\frac25 v(B)\Big\rc.$$
Thus $H'$ is a simple good graph. Since $H\setm G_0$ is very good, the assumption of Lemma \ref{glue} is satisfied by $H\setm G_0$  and $H'$. 
Based on this lemma, after a few further rounds the host graph is a good graph $G$ with 
$$
\ess(G)=\ess(H\setm G_0)+\ess(H')=\ess(H)+2=m+r=n.
$$ 
Further analysis is the same as in Stage 3 
for non-exceptional cases above. It leads to the conclusion that a blue path of length $n-1$ is created
and the number of edges on the board at this moment is not greater than $\lc 7n/5\rc-1$.
 
\smallskip\noindent
\textbf{Subcase 2.} $G_0=\emptyset$ and $G_1\neq\emptyset$.\smallskip

Let  $u_1,u'_1$ be the 1-ends of $G_1$.
Builder forces blue edges $u_1a$ and $u'_1b$. 
As a result we obtain a blue $(1,1)$-path $G'_1$ with 1-ends $a$ and $b$, on $v(G_1)+v(P)$ vertices.
Let $G=H\cup P\cup\{u_1a,u'_1b\}$. 
Since $5\mid v(G_1\cup G_2)$, the subgraph $G'_1\cup G_2$ of the coloured graph $G$ satisfies 
\begin{eqnarray*}
|\rn_{G}(G'_1\cup G_2)|&\le& |\rn_H(G_1\cup G_2)|+|\rn_G(P)|\le \frac25 v(G_1\cup G_2)+1=
\frac25 v(G_1\cup G_2)+\Big\lc\frac25 v(P)\Big\rc\\
&=&\Big\lc\frac25(v(G_1\cup G_2)+v(P))\Big\rc=\Big\lc\frac25(v(G'_1\cup G_2))\Big\rc.
\end{eqnarray*}
It is not hard to verify that $G$ satisfies all conditions of a good graph
and its essential subgraphs are $(\emptyset, G'_1,G_2,L)$.  
Again we use the same argument as in Stage 3 
for non-exceptional cases and conclude the assertion. 

\smallskip\noindent
\textbf{Subcase 3.} $G_0=\emptyset$ and $G_1=\emptyset$.\smallskip
 
Then $G_2\neq\emptyset$ since $G$ has at least 10 essential vertices. 
 Let $w_2$ be the transition vertex of $G_2$ and $u_2$ be its red end.
 Let us recall that the edge $w_2u_2$ is red.  
Builder forces blue edges: $w_2a$, $au_2$ and $u_2b$. 
Thereby we obtain a blue $(1,0)$-path $F$ on $v(G_2)+v(P)$ vertices and its 1-end $b$. 

Apart from the extended $(2,0)$-path $G'_2$, $H$ may contain also a limb.  Then, as in 
Stage 3 for non-exceptional case,
Builder transforms the limb into a blue $(1,0)$-path $P'$ with an 1-end $x_4$ and 
$|\rn(P')|\le 2=\frac25 v(L)$. 
Afterwards, Builder forces the edge $x_4b$ blue and a longer blue path $B$ is obtained.
 
Let $B=F\cup P'\cup\{x_4a\}$ if $L\neq\emptyset$, while otherwise we put $B=F$ and define $P'=\emptyset$.
Note that $|\rn(P')|\le\frac25 v(L)$  in both cases ($L=\emptyset$ or not).  
Notice also that in both cases $V(B)=V(G_2\cup L \cup P)$ so $v(B)=\ess(H)+2=n+r=n$. So there is a blue path
on $n$ vertices on the board and Builder wins.
  
Furthermore, since $5\mid v(G_2\cup G_1)=v(G_2)$, we have  
\begin{eqnarray*}
|\rn(B)|&\le& |\rn_H(G_2)|+|\rn(P)|+|\rn(P')|\le \frac25 v(G_2)+1+\frac25 v(L)\\
&=&\frac25 v(G_2)+\Big\lc\frac25 v(P)\Big\rc+\frac25 v(L)
=\Big\lc\frac25 v(B)\Big\rc.
\end{eqnarray*}
It follows from the strategy of Builder that $\rn(B)$ consists of all red edges of the game $\RR_\emptyset(P_4,P_n)$ and the blue path $B$ contains all blue edges present on the board. 
Thus the number of rounds in the game is 
$$e(B)+|\rn(B)|\le \Big\lc\frac75 v(B)\Big\rc-1=\Big\lc\frac75 n\Big\rc-1.$$
This ends the proof in the exceptional case and the proof of Theorem \ref{mainp4}. 

\bibliographystyle{amsplain}

\end{document}